\documentclass[3pt]{article}

\usepackage{graphicx}\usepackage{epstopdf}\usepackage{color}
\usepackage{amsmath}\usepackage{amsfonts}\usepackage{amssymb}
 \usepackage{amsthm}

\numberwithin{equation}{section}

\begin{document}


\title{A New Multigrid Finite Element Method  for the Transmission Eigenvalue Problems}

\author{ {Jiayu Han, Yidu Yang, Hai Bi} \\\\
{\small School of Mathematics and Computer Science, }\\{\small
Guizhou Normal University,  Guiyang,  $550001$,  China}\\{\small
hanjiayu126@126.com, ydyang@gznu.edu.cn, bihaimath@gznu.edu.cn}}
\date{~}
\pagestyle{plain} \textwidth 145mm \textheight 215mm \topmargin 0pt
\maketitle

\indent{\bf\small Abstract~:~} {\small Numerical methods for the transmission
eigenvalue problems are hot topics  in recent years.
Based on the work of  Lin and Xie [Math.   Comp., 84(2015), pp. 71-88], we build a  multigrid method  to solve the  problems.   With our method, we only need to solve a series of primal and dual eigenvalue problems on a coarse mesh and the associated boundary value problems
on the finer and finer meshes. Theoretical analysis and   numerical results show that our method is simple and easy to implement and is efficient for computing real and complex transmission eigenvalues.
}\\
\indent{\bf\small Keywords~:~\scriptsize} {\small
Transmission eigenvalues, Multigrid method, Nonsymmetric eigenvalue problems, Extended/Generalized finite element.}
\section{Introduction}
\label{intro}
\indent
The  transmission eigenvalue problems have theoretical importance in the uniqueness and reconstruction in
inverse scattering theory \cite{cakoni,colton1}.  Transmission eigenvalues
can be determined from the far-field data of the scattered wave and used to obtain estimates
for the material properties of the scattering object \cite{cakoni1,cakoni2}. Many literatures such as
\cite{colton1,cakoni2,kirsch,paivarinta,rynne} studied
the existence of transmission eigenvalues, and \cite{cakoni2,cakoni3,colton3} et al.
explored the upper and lower bounds for the index of refraction
$n(x)$.

 In recent years, numerical methods for the transmission eigenvalue problems have  attracted the attention from more and more researchers.
 The first numerical treatment of the transmission eigenvalue problem appeared in
\cite{colton2} where three finite element methods are proposed for the Helmholtz transmission
eigenvalues. Later on, many other numerical methods were developed to solve the problems  (see, e.g., 
 \cite{sun1,ji2,an1,cakoni4,yang1}).  In particular,  Sun \cite{sun1} proposed an iterative method  and gave a coarse  error analysis. Furthermore, Ji et al. \cite{ji2} developed his work and  proved the accurate error estimates for both eigenvalues and eigenfunctions  by constructing an auxiliary problem as a bridge. An and Shen \cite{an1} proposed a spectral-element method to solve this problem numerically.
 Afterwards using the linearized technique the authors in \cite{cakoni4,yang1}   builded two new weak formulations and the corresponding finite element discretizations.

The idea of multigrid methods for eigenvalue problems was developed originally from two grid methodology.  In 2001, Xu and Zhou \cite{xu1} proposed
 a two grid method based on inverse iteration for elliptic eigenvalue problems, which is, in a way, related to that in
Lin and Xie \cite{lin3}. After that, the two grid method was further developed into multigrid method and local and parallel algorithm  in \cite{xu2,yang2,bi,dai} et al.
  In recent years, Lin and Xie \cite{lin,xie} proposed
a mutilevel correction method. This method can be regarded as  the combination of  two grid method  and the extended/generalized finite element method.
  The extended/generalized finite element method   was developed in 1990s  by   \cite{moes,belytschko1,duarte,melenk} et al, which has important applications on problems  in material science \cite{belytschko,fries}.
The method of Lin and Xie \cite{lin,xie} enriches the finite element space at each correction step with  the numerical eigenfunctions obtained from the last  step.
So it is able to naturally reproduce the   feature of eigenfunctions: the discontinuity, singularity, boundary layer, etc. Such an embedding of the problem's feature into the finite element space can significantly improve convergence rates and accuracy step by step, so that
 the multigrid method can achieve the same accuracy as solving the eigenvalue problem directly but with less computational work.

In this paper, based on the literatures \cite{lin,xie}, we propose a new  multigrid method  to solve the transmission eigenvalue problems 
but based on the new weak formulation (\ref{s2.6}) proposed in \cite{yang1} which is a linear and nonsymmetric eigenvalue problem.  In this work, (1) we prove the error estimates of the  transmission eigenvalues and eigenfunctions for our multigrid method. Our theoretical results 
 are valid for arbitrary real and complex eigenvalues.
 (2) 
 With our method, due to adopting the linearized weak formulation, we can transform the transmission eigenvalue problem  into a generalized matrix eigenvalue problem and can be solved efficiently by the sparse solver $eigs$ in Matlab;
 (3) with our multigrid method, the solution of eigenvalue problem on a fine mesh can be reduced to  a series of the solutions of the eigenvalue problem on a coarse meshes and a series of solutions of the boundary value problems on the multilevel meshes. 
 As    numerical results indicate, this method is applicable to the real and complex transmission eigenvalues.

\section{Preliminaries}
\indent Let $H^{s}(D)$ be a Sobolev space with norm $\|\cdot\|_{s}$ ($s=1,2$),
  and  $$H_{0}^{2}(D)=\{v\in H^{2}(D):
v|_{\partial D}=\frac{\partial v}{\partial \nu}|_{\partial D}=0\}.$$

 \indent Consider the Helmholtz transmission eigenvalue
problem: Find $k\in \mathcal{C}$, $\omega, \sigma\in L^{2}(D)$,
$\omega-\sigma\in H^{2}(D)$ such that
\begin{eqnarray}\label{s2.1}
&&\Delta \omega+k^{2}n(x)\omega=0,~~~in~ D,\\\label{s2.2}
 &&\Delta
\sigma+k^{2}\sigma=0,~~~in~ D,\\\label{s2.3}
 &&\omega-\sigma=0,~~~on~ \partial
D,\\\label{s2.4}
 &&\frac{\partial \omega}{\partial \nu}-\frac{\partial
\sigma}{\partial \nu}=0,~~~ on~\partial D,
\end{eqnarray}
where $D\subset \mathbb{R}^{2}$ or $D\subset \mathbb{R}^{3}$  is a bounded simply
connected inhomogeneous medium,
 $\nu$ is the unit outward normal to $\partial D$.\\

It is possible to write (\ref{s2.1})-(\ref{s2.4})
as an equivalent eigenvalue problem for  $u=\omega-\sigma\in H_{0}^{2}(D)$. In particular,
\begin{eqnarray*}
 (\Delta u +k^2u)=\Delta \omega+k^2 \omega=k^2 (1-n)\omega.
\end{eqnarray*}
Dividing by $n-1$ and applying the operator $(\Delta+k^2n)$ to the above equality,
 the
eigenvalue problem (\ref{s2.1})-(\ref{s2.4}) can be stated as
follows: Find $k^{2}\in \mathcal{C}$, $k^{2}\not=0$, nontrivial $u\in
H_{0}^{2}(D)$ such that
\begin{eqnarray}\label{s2.5}
(\frac{1}{n(x)-1}(\Delta u+k^{2} u),\Delta v+k^{2} n(x)v
)_{0}=0,~~~\forall v \in~H_{0}^{2}(D),
\end{eqnarray}
where $(\cdot,\cdot)_0$ is the inner product of $L^2(D)$.
  As usual, we  define $\lambda=k^{2}$ as   the transmission eigenvalue in this paper.
We suppose that the index of refraction $n\in
L^{\infty}(D)$ satisfying either one of the following assumptions
\begin{eqnarray*}
&&(C1)~~1+\delta\leq \inf_{D}n(x)\leq n(x)\leq
\sup_{D}n(x)<\infty,\\
&&(C2)~~0< \inf_{D}n(x)\leq n(x)\leq \sup_{D}n(x)<1-\varrho,
\end{eqnarray*}
for some constant $\delta>0$ or $\varrho>0$.\\

\indent For simplicity, in the coming discussion we assume  $(C1)$
holds and $n(x)$ is proper smooth (for example $n(x)\in W^{2,\infty}(D)$). For the case $(C2)$  the argument method is the same.

\indent Define Hilbert space $\mathbf{H}=H_{0}^{2}(D)\times
L^{2}(D)$ with norm
$\|(u,w)\|_{\mathbf{H}}=\|u\|_{2}+\|w\|_{0}$, and define
$\mathbf{H}_{s}=H^{s}(D)\times H^{s-2}(D)$ with norm
$\|(u,w)\|_{\mathbf{H}_{s}}=\|u\|_{s}+\|w\|_{s-2}$, $s=0,1$.\\

\indent From (\ref{s2.5}) we derive that
\begin{eqnarray*}
&&(\frac{1}{n-1}\Delta u, \Delta v)_{0}-\lambda
(\nabla(\frac{1}{n-1}u), \nabla v)_{0}\nonumber\\
&&~~~~~~-\lambda(\nabla u, \nabla(\frac{n}{n-1}v))_{0}+\lambda^{2}
(\frac{n}{n-1}u,v )_{0}=0,~~~\forall v \in~H_{0}^{2}(D).
\end{eqnarray*}
 Let $w=\lambda u$,
we arrive at a linear
weak formulation: Find $(\lambda, u, w)\in \mathcal{C}\times
H_{0}^{2}(D)\times L^{2}(D)$ such that
\begin{eqnarray*}
&&(\frac{1}{n-1}\Delta u, \Delta v)_{0}=\lambda
(\nabla(\frac{1}{n-1}u),
\nabla v)_{0}\nonumber\\
&&~~~~~~+\lambda(\nabla u, \nabla(\frac{n}{n-1}v))_{0}-\lambda
(\frac{n}{n-1}w, v )_{0},~~~\forall v
\in~H_{0}^{2}(D),\\
 &&(w,
z)_{0}=\lambda(u,z)_{0},~~~\forall z\in L^{2}(D).
\end{eqnarray*}

\indent We introduce the following sesquilinear forms
\begin{eqnarray*}
&&A((u,w),(v,z))=(\frac{1}{n-1}\Delta u, \Delta v)_{0}+(w, z)_{0},\\
&&B((u,w),(v,z))\nonumber\\
&&~~~=(\nabla(\frac{1}{n-1}u), \nabla v)_{0}+(\nabla u,
\nabla(\frac{n}{n-1}v))_{0}- (\frac{n}{n-1}w, v
)_{0}+(u,z)_{0}\\
&&~~~=-(\frac{1}{n-1}u, \nabla\cdot\nabla v)_{0}-(u,
\nabla\cdot\nabla(\frac{n}{n-1}v))_{0}- (\frac{n}{n-1}w, v
)_{0}+(u,z)_{0},
\end{eqnarray*}
then (\ref{s2.5})  can be rewritten as: Find
$\lambda\in \mathcal{C}$,  nontrivial $(u,w)\in \mathbf{H}$ such that
\begin{eqnarray}\label{s2.6}
A((u,w),(v,z)) =\lambda B((u,w),(v,z)),~~~\forall (v,z)\in
\mathbf{H}.
\end{eqnarray}
Let norm $\|\cdot\|_A $ be induced by the inner product $A(\cdot,\cdot)$, then it is clear $\|\cdot\|_A$ is equivalent to $\|\cdot\|_\mathbf{H}$. \\




\indent One can easily verify that for any
given $(f,g)\in \mathbf{H}_{s}$ ($s=0,1$), $B((f,g), (v, z))$ is a
continuous linear form on
$\mathbf{H}$: $$B((f,g), (v, z))\lesssim \|(f,g)\|_{\mathbf{H}_s}\|(v,z)\|_{\mathbf{H}},~\forall (v,z)\in \mathbf{H}.$$
Here and hereafter this paper, we use the symbols  $x \lesssim y$ to mean $x \le Cy$ for a constant $C$ that is independent of mesh size and iteration times and may be different at different occurrences.\\
\indent The source problem associated with (\ref{s2.6}) is given
by: Find $(\psi,\varphi)\in \mathbf{H}$ such that
\begin{eqnarray}\label{s2.7}
A((\psi,\varphi),(v,z)) =B((f,g),(v,z)),~~~\forall (v,z)\in
\mathbf{H}.
\end{eqnarray}
From the Lax-Milgram theorem we know that the problem (\ref{s2.7}) exists an
unique solution, therefore, we define the corresponding solution
operator  $T: \mathbf{H}_{s}\to \mathbf{H}$ by
\begin{eqnarray}\label{s2.8}
A(T(f,g),(v,z)) =B((f,g), (v,z)),~~~\forall (v,z)\in \mathbf{H}.
\end{eqnarray}
Then (\ref{s2.6}) has the equivalent operator form:
\begin{eqnarray}\label{s2.9}
T(u,w)= {\lambda}^{-1}(u,w).
\end{eqnarray}
\indent Consider the dual problem of (\ref{s2.6}): Find
$\lambda^{*}\in \mathcal{C}$, nontrivial $(u^{*},w^{*})\in \mathbf{H}$
such that
\begin{eqnarray}\label{s2.10}
A((v,z), (u^{*},w^{*})) =\overline{\lambda^{*}} B((v,z),
(u^{*},w^{*})),~~~\forall (v,z)\in \mathbf{H}.
\end{eqnarray}

Note that the primal and dual eigenvalues are connected via
$\lambda=\overline{\lambda^{*}}$.\\

Define the corresponding solution operator  $T^{*}:
\mathbf{H}_{s}\to \mathbf{H}$ by
\begin{eqnarray}\label{s2.11}
A((v,z), T^{*}(f,g)) =B((v,z), (f,g)),~~~\forall (v,z)\in
\mathbf{H}.
\end{eqnarray}
Then (\ref{s2.10}) has the equivalent operator form:
\begin{eqnarray}\label{s2.12}
T^{*}(u^{*},w^{*})=\lambda^{*-1}(u^{*},w^{*}).
\end{eqnarray}

\indent Clearly,   $T^{*}$ is the adjoint operator of $T$
in the sense of inner product $A(\cdot,\cdot)$. 
In order to discretize the space
$\mathbf{H}$, we need  two finite element spaces to  discretize $ H_{0}^{2}(D)$ and
$L^{2}(D)$, respectively, but here we can construct only one conforming finite element space $S^h\subset H_{0}^{2}(D)$ such that $\mathbf{H}_{h}:=S^h\times S^h\subset H_{0}^{2}(D) \times L^{2}(D)$.

\indent The conforming finite element approximation of (\ref{s2.6})
is given by the following: Find $\lambda_{h}\in \mathcal{C}$, nontrivial
$(u_{h},w_{h})\in \mathbf{H}_{h}$ such that
\begin{eqnarray}\label{s2.13}
A((u_{h},w_{h}),(v,z)) =\lambda_{h}
B((u_{h},w_{h}),(v,z)),~~~\forall (v,z)\in \mathbf{H}_{h}.
\end{eqnarray}

We introduce the corresponding solution operator: $T_{h}:
\mathbf{H}_{s}\to \mathbf{H}_{h}$ (s=0,1):
\begin{eqnarray}\label{s2.15}
A(T_{h}(f,g),(v,z)) =B((f,g), (v,z)),~~~\forall (v,z)\in
\mathbf{H}_{h}.
\end{eqnarray}
Then (\ref{s2.13}) has the operator form:
\begin{eqnarray}\label{s2.16}
T_{h}(u_{h},w_{h})={\lambda_{h}^{-1}}(u_{h},w_{h}).
\end{eqnarray}

\indent Define the projection operators $P_{h}^{1}: H_{0}^{2}(D)\to
S^{h}$ and $P_{h}^{2}: L^{2}(D)\to S^{h}$ by
\begin{eqnarray}\label{s2.17}
&&(\frac{1}{n-1}\Delta (u-P_{h}^{1}u), \Delta v)_{0}=0,~~~\forall
v\in S^{h},\\\label{s2.18} &&(w-P_{h}^{2}w,
z)_{0}=0,~~~\forall z\in S^{h}.
\end{eqnarray}
Let
$$P_{h}(u,w)=(P_{h}^{1}u, P_{h}^{2}w),~~~\forall (u,w)\in \mathbf{H}.$$
Then $P_{h}: \mathbf{H}\to \mathbf{H}_{h}$, and
\begin{eqnarray}\label{s2.19}
&&A((u,w)-P_{h}(u,w),(v,z))
=0,~~~\forall (v,z)\in \mathbf{H}_{h}.
\end{eqnarray}

 \indent We need the following regularity assumption:\\
\indent {\bf R(D)}.~~ For any $\xi\in H^{-s}(D)$ ($s=0,1$), there
exists $\psi\in H^{2+r_{s}}(D)$ satisfying
\begin{eqnarray*}
&&\Delta(\frac{1}{n-1}\Delta \psi)=\xi,~~~in~D,\\
&&\psi=\frac{\partial \psi}{\partial \nu}=0~~~ on~\partial D,
\end{eqnarray*}
and
\begin{eqnarray}\label{2.19}
\|\psi\|_{2+r_{s}}\leq C_{p} \|\xi\|_{-s},~~~s=0,1
\end{eqnarray}
where $r_{1}\in (0,1]$, $r_{0}\in (0,2]$, $C_{p}$ denotes the prior
constant dependent on the equation and $D$ but
independent of the right-hand side $\xi$ of the equation.

It is easy to know that  (\ref{2.19})  is valid with $r_s = 2-s$ when $n \in W^{2,p}(D)$
($p$ is greater than but arbitrarily close to 2) and $\partial D$ is appropriately smooth.
When $D \subset R^2$ is a convex polygon, from Theorem 2 in \cite{blum}, when $n \in W^{2,p}(D)$,
we can get $r_1 = 1$ and that if the inner angle at each critical boundary point is
smaller than 126.283696...0 then $r_0 = 2$.

The following Lemmas 2.1-2.3 come from \cite{yang1}. They give   in the sense of   lower norms the estimates of the finite element projection and the convergence of $T_h$ to $T$.\\
\indent{\bf Lemma 2.1 \cite[Lemma 3.4]{yang1}.~~}Suppose that  $n\in W^{2,\infty}(D)$ and \textbf{R(D)}   is valid
$(s=0,1)$, then
 for $(u, w)\in \mathbf{H}$,
\begin{eqnarray}\label{2.20}
\|(u, w)-P_{h}(u,w)\|_{\mathbf{H}_{s}} \lesssim
h^{r_{s}}\| (u,w)-P_{h}(u,w)\|_{\mathbf{H}},~~~s=0,1.
\end{eqnarray}

\indent The conforming finite element approximation of (\ref{s2.10})
is given by: Find $\lambda_{h}^{*} \in \mathcal{C}$,
$(u_{h}^{*},w_{h}^{*})\in \mathbf{H}_{h}$ such that
\begin{eqnarray}\label{s2.20}
A((v,z), (u_{h}^{*},w_{h}^{*})) =\overline{\lambda_{h}^{*}}
B((v,z), (u_{h}^{*},w_{h}^{*})),~~~\forall (v,z)\in
\mathbf{H}_{h}.
\end{eqnarray}

Note that the primal and dual eigenvalues are connected via
$\lambda_{h}=\overline{\lambda_{h}^{*}}$.\\

Define the solution operator $T_{h}^{*}: \mathbf{H}_{s}\to
\mathbf{H}_{h}$ satisfying
\begin{eqnarray}\label{s2.21}
A((v,z), T_{h}^{*}(f,g))=B((v,z), (f,g)),~~~\forall~(v,z)\in
\mathbf{H}_{h}.
\end{eqnarray}
Naturally (\ref{s2.20}) has the following equivalent operator form
\begin{eqnarray}\label{s2.22}
T_{h}^{*}(u_{h}^{*},w_{h}^{*})=\lambda_{h}^{*-1}(u_{h}^{*},w_{h}^{*}).
\end{eqnarray}

\indent{\bf Lemma 2.2 \cite[Theorem 3.1]{yang1}.~~}Let $n \in W^{1,\infty}(D)$, then
\begin{eqnarray}\label{s2.23}
&&\|T-T_{h}\|_{\mathbf{H}}\to 0,\\\label{s2.24}
&&\|T-T_{h}\|_{\mathbf{H}_{1}}\to 0,
\end{eqnarray}
and let $n\in W^{2,\infty}(D)$, then
\begin{eqnarray}\label{s2.25}
\|T-T_{h}\|_{\mathbf{H}_{0}}\to 0.
\end{eqnarray}

In this paper, let $\lambda_i$ be the $ith$ eigenvalue of (\ref{s2.6}) with the
algebraic multiplicity $q$ and the ascent  1. Then, according to spectral approximation theory \cite{babuska1,chatelin}, there are  $q$
eigenvalues $\lambda_{j,h}$ ($j=i,\cdots,i+q-1$)  of (\ref{s2.13})   converging to $\lambda_i$. Let $M(\lambda_i)$ be the
space spanned by all eigenfunctions corresponding to the
eigenvalue $\lambda_i$. Let $M_h(\lambda_i)$ be the space spanned by all
generalized eigenfunctions corresponding to the numerical eigenvalues $\{\lambda_{j,h}\}_{j=i}^{i+q-1}$ of (\ref{s2.13}). As for the dual problems (\ref{s2.10}) and (\ref{s2.20}), the definitions
of $M^*(\lambda^*_i)$  and $M^*_h(\lambda^*_i)$
  are made similarly to $M(\lambda_i)$  and $M_h(\lambda_i)$, respectively.

In what follows, to describe the approximation  relation between the finite element space $\mathbf{H}_h$ and the eigenfunction spaces $ M({\lambda}_i)$ and $M^*({\lambda}^*_i)$, we introduce the following quantities
 \begin{eqnarray*}
 && {\delta}_{h}({{\lambda}}_i)=\sup\limits_{(v,z)\in M({\lambda}_i)\atop\|(v,z)\|_\mathbf{H}=1}\inf\limits_{(v_h,z_h)\in \mathbf{H}_{h}}\|(v,z)-(v_h,z_h)\|_\mathbf{H} , \\
 && {\delta}^*_{h}({{\lambda}}^*_i)=\sup\limits_{(v,z)\in M^*({\lambda}^*_i)\atop\|(v,z)\|_\mathbf{H}=1}\inf\limits_{(v_h,z_h)\in \mathbf{H}_{h}}\|(v,z)-(v_h,z_h)\|_\mathbf{H} .
 \end{eqnarray*}


The operator convergence results  (\ref{s2.23})-(\ref{s2.25}) are critical as  a bridge of making the error analysis for the discrete problem (\ref{s2.13}). From these results, we yield immediately the following lemma using the spectral approximation theory. \\

\indent{\bf Lemma 2.3.~~}  Suppose $n \in W^{2,\infty}(D)$.
Let $(u_{j,h},w_{j,h})$ ($j=i,i+1,\cdots,i+q-1$) be eigenfunction corresponding to
$\lambda_{j,h}$ and $\|(u_{j,h},w_{j,h})\|_A=1$, then there exists eigenfunction
$(u_i,w_i)$ corresponding to $\lambda_i$ such that
\begin{eqnarray}\label{s2.28}
&&\|(u_{j,h},w_{j,h})-(u,w)\|_{\mathbf{H}}\lesssim
\delta_{h}(\lambda_i) ,\\\label{s2.29}
&&\|(u_{j,h},w_{j,h})-(u,w)\|_{\mathbf{H}_s}\lesssim
h^{r_s}\delta_{h}(\lambda_i) ,~s=0,1,\\\label{s2.30}
&&|\lambda_i-\lambda_{j,h}|\lesssim
 \delta_{h}(\lambda_i)\delta_{h}^{*}(\lambda_i^{*}).
\end{eqnarray}
\indent{\bf Remark 2.1.}~~The similar estimates as above are valid for the dual problem (\ref{s2.13}) (see \cite{yang1}).

\section{A New Multigrid Method}
\label{sec:2}
In this section, based on the multilevel correction method proposed by Lin and Xie \cite{lin,xie}, we give the multigrid scheme for the weak form (\ref{s2.6}).
 Our theoretical results are given in Theorems 3.1-3.2.  Prior to our argument, we give the following basic condition related to finite element spaces and their    approximation relation to eigenfunction spaces.

%
We construct the finite element spaces such that
$\mathbf{H}_{H}=\mathbf{H}_{h_{1}}\subset\mathbf{H}_{h_{2}}\subset\cdots\subset\mathbf{H}_{h_{n}}$
and
\begin{eqnarray}\label{s3.1}
\delta_{h_{m+1}}\approx\frac{1}{\beta}\delta_{h_{m}}, ~\delta^*_{h_{m+1}}\approx\frac{1}{\beta}\delta^*_{h_{m}},
\end{eqnarray}
where $\beta>1$ is a constant only dependent on the smoothness of eigenfunctions corresponding to $\lambda_i$  and the degree $t$ of the piecewise polynomial space.

   The above condition is readily satisfied  on   regular meshes. In particular, 
  if the meshes are obtained from a procedure of bisection mesh refinement and $M(\lambda_i),
M^*(\lambda^*_i)\subset H^{2+r}(D)\times H^{r}(D)$ ($r\le t-1$),   then we have approximately  $\beta\approx 2^r$.

Assume that we have obtained the eigenpair approximations $(\lambda_{j,h_m},u_{j,h_{m}},w_{j,h_{m}})$ $(j=i,i+1,\cdots,i+q-1)$ and the corresponding dual ones
$(\lambda^*_{j,h_m},u^*_{j,h_{m}},w^*_{j,h_{m}})$. First of all, we give  one correction step of the multigrid scheme.\\

\indent{\bf Algorithm 1.~~}One Correction Step.\\
\textbf{Step 1}. Solve the following linear boundary value problems: For $j=i,i+1,\cdots,i+q-1$, find $(\widehat{u}_{j,h_{m+1}},\widehat{w}_{j,h_{m+1}})\in \mathbf{H}_{h_{m+1}}$ such
that
\begin{eqnarray*}
A((\widehat{u}_{j,h_{m+1}},\widehat{w}_{j,h_{m+1}}),
(v,z))=\lambda_{j,h_m}B((u_{j,h_{m}},w_{j,h_{m}}),(v,z)),~~~\forall (v,z)\in
\mathbf{H}_{h_{m+1}},
\end{eqnarray*}
and find $(\widehat{u}^*_{j,h_{m+1}},\widehat{w}^*_{j,h_{m+1}})\in \mathbf{H}_{h_{m+1}}$ such that
\begin{eqnarray*}
A(
(v,z),(\widehat{u}^*_{j,h_{m+1}},\widehat{w}^*_{j,h_{m+1}}))=\lambda_{j,h_m}
B((v,z),(u^*_{j,h_{m}},w^*_{j,h_{m}})),~~~\forall (v,z)\in
\mathbf{H}_{h_{m+1}}.
\end{eqnarray*}
\textbf{Step 2}. Construct a new finite element space
$\mathbf{H}_{H,h_{m+1}}\supseteq\mathbf{H}_{H}+span\Big\{(\widehat{u}_{j,h_{m+1}},\widehat{w}_{j,h_{m+1}}),$ $(\widehat{u}^*_{j,h_{m+1}},\widehat{w}^*_{j,h_{m+1}}),$ $j=i,i+1,\cdots,i+q-1\Big\}$ and solve the following eigenvalue problems:
For $j=i,i+1,\cdots,i+q-1$,
find $ \lambda_{j,h_{m+1}}\in \mathcal{C},$ nontrivial $( {u}_{j,h_{m+1}}, {w}_{j,h_{m+1}}) \in \mathbf{H}_{H,h_{m+1}}$ such
that
\begin{eqnarray*}
A(( {u}_{j,h_{m+1}}, {w}_{j,h_{m+1}}),
(v,z))=\lambda_{j,h_{m+1}}B((u_{j,h_{m}},w_{j,h_{m+1}}),(v,z)),~~   \forall (v,z)\in
\mathbf{H}_{H,h_{m+1}},
\end{eqnarray*}
and find nontrivial $ ( {u}^*_{j,h_{m+1}}, {w}^*_{j,h_{m+1}})\in \mathbf{H}_{H,h_{m+1}}$ such that
\begin{eqnarray*}
A(
(v,z),( {u}^*_{j,h_{m+1}}, {w}^*_{j,h_{m+1}}))=\lambda_{j,h_{m+1}}
B((v,z),(u^*_{j,h_{m+1}},w^*_{j,h_{m+1}})),~~\forall (v,z)\in
\mathbf{H}_{H,h_{m+1}}.
\end{eqnarray*}
\indent We output $\Big\{\lambda_{j,h_{m+1}}\Big\}_{j=i}^{i+q-1}$
  and  a basis $\Big\{(u_{j,h_{m+1}},w_{j,h_{m+1}})\Big\}_{j=i}^{i+q-1}$ of
$M_{h_{m+1}}(\lambda_i)$ with $\|(u_{j,h_{m+1}},w_{j,h_{m+1}})\|_ {A}=1$ and a basis $\Big\{(u^*_{j,h_{m+1}},w^*_{j,h_{m+1}})\Big\}_{j=i}^{i+q-1}$ of
$M^*_{h_{m+1}}(\lambda^*_i)$  with $\|(u^*_{j,h_{m+1}},w^*_{j,h_{m+1}})\|_ {A}=1$.
We define the above two steps as Procedure $\mathbf{Correction}$:
\begin{eqnarray*}
&&\Big\{\lambda_{j,h_{m+1}},u_{j,h_{m+1}},w_{j,h_{m+1}},u^*_{j,h_{m+1}},w^*_{j,h_{m+1}}\Big\}_{j=i}^{i+q-1}\\
&&=Correction\Big(\mathbf{H}_H,\Big\{\lambda_{j,h_{m}},u_{j,h_{m}},w_{j,h_{m}},u^*_{j,h_{m}},w^*_{j,h_{m}}\Big\}_{j=i}^{i+q-1},\mathbf{H}_{h_{m+1}}\Big)
\end{eqnarray*}

Implementing Procedure  $\mathbf{Correction}$ repeatedly leads to the following Scheme.
\\\\
\textbf{Algorithm 2.~~}Mutigrid Scheme.\\
\textbf{Step 1}. Construct a series of nested finite element spaces $\mathbf{H}_H=\mathbf{H}_{h_1},$ $\mathbf{H}_{h_2},$ $\cdots,$ $\mathbf{H}_{h_n}$
such that (\ref{s3.1}) holds.\\
\textbf{Step 2}. For $j=i,i+1,\cdots,i+q-1$,
find $(\lambda_{j,h_{1}}, {u}_{j,h_{1}}, {w}_{j,h_{1}}) \in \mathcal{C}\times\mathbf{H}_{h_{1}}$ such
that $\|({u}_{j,h_{1}}, {w}_{j,h_{1}})\|_A=1$ and
\begin{eqnarray*}
A(( {u}_{j,h_{1}}, {w}_{j,h_{1}}),
(v,z))=\lambda_{j,h_{1}}B((u_{j,h_{1}},w_{j,h_{1}}),(v,z)),~~~\forall (v,z)\in
\mathbf{H}_{h_{1}},
\end{eqnarray*}
and find $({u}^*_{j,h_{1}}, {w}^*_{j,h_{1}}) \in  \mathbf{H}_{h_{1}}$ such that $\|({u}^*_{j,h_{1}}, {w}^*_{j,h_{1}})\|_A=1$ and
\begin{eqnarray*}
A(
(v,z),( {u}^*_{j,h_{1}}, {w}^*_{j,h_{1}}))=\lambda_{j,h_{1}}
B((v,z),(u^*_{j,h_{1}},w^*_{j,h_{1}})),~~~\forall (v,z)\in
\mathbf{H}_{h_{1}}.
\end{eqnarray*}
\textbf{Step 3}. For $m=1,2,\cdots,N-1$\\
Obtain a new eigenpair approximation $(\lambda_{j,h_{m+1}},u_{j,h_{m+1}},w_{j,h_{m+1}})$ by
\begin{eqnarray*}
&&\Big\{\lambda_{j,h_{m+1}},u_{j,h_{m+1}},w_{j,h_{m+1}},u^*_{j,h_{m+1}},w^*_{j,h_{m+1}}\Big\}_{j=i}^{i+q-1}\\
&&=Correction\Big(\mathbf{H}_H,\Big\{\lambda_{j,h_{m}},u_{j,h_{m}},w_{j,h_{m}},u^*_{j,h_{m}},w^*_{j,h_{m}}\Big\}_{j=i}^{i+q-1},\mathbf{H}_{h_{m+1}}\Big)
\end{eqnarray*}
End
\\

Before making the error analysis of multigrid scheme, we give the following assumption.
\\

%


{\indent{\bf (A) Quasi-biorthogonality}.~~
Suppose that there are $\Big\{(\widetilde{u}_{j,h_{m}},\widetilde{w}_{j,h_{m}})\Big\}_{j=i}^{i+q-1}$ $\subset M_{h_m}(\lambda_i)$ with $\|(\widetilde{u}_{j,h_{m}}$, $\widetilde{w}_{j,h_{m}})\|_\mathbf{H}$ $=1$ and $\Big\{(\widetilde{u}^*_{j,h_{m}},\widetilde{w}^*_{j,h_{m}})\Big\}_{j=i}^{i+q-1}\in M^*_{h_m}(\lambda^*_i)$ with $\|(\widetilde{u}^*_{j,h_{m}}$, $\widetilde{w}^*_{j,h_{m}})\|_\mathbf{H}$ $=1$ such that $$|A(({u}_{j,h_{m}},{w}_{j,h_{m}}),(\widetilde{u}^*_{l,h_{m}},\widetilde{w}^*_{l,h_{m}}))|+|A((\widetilde{u}_{j,h_{m}},\widetilde{w}_{j,h_{m}}),({u}^*_{l,h_{m}},{w}^*_{l,h_{m}}))|\lesssim H^{r_0},$$ ($j,l=i,i+1,\cdots,i+q-1,j\neq l$),\\
and $|A(({u}_{j,h_{m}},{w}_{j,h_{m}}),(\widetilde{u}^*_{j,h_{m}},\widetilde{w}^*_{j,h_{m}}))|+|A((\widetilde{u}_{j,h_{m}},\widetilde{w}_{j,h_{m}}),({u}^*_{j,h_{m}},{w}^*_{j,h_{m}}))|$
 $(j=i,i+1,\cdots,i+q-1)$ has a positive lower bound uniformly with respective to $h_m$.
\\\\
\indent When $\lambda_i$ is a simple eigenvalue ($q= 1$), it is clear that Assumption ($\mathbf{A}$) is valid; when $q>1$, we can prove the follow conclusion:\\
\indent If the distance (in $\|\cdot\|_\mathbf{H}$) from $({u}_{j,h_{m}},{w}_{j,h_{m}})$ ($j=i,i+1,\cdots,i+q-1$) to $span\Big\{({u}_{l,h_{m}},{w}_{l,h_{m}}),$ $l=i,i+1,\cdots,i+q-1,l\neq j\Big\}$ has a positive lower bound uniformly with respective to $h_m$,
then Assumption ($\mathbf{A}$) is valid.\\\\
\indent The following theorem, an  extension of the corresponding theorems in \cite{lin,xie,peng}, indicates that the accuracy of numerical eigenpair can be apparently improved after one correction step.\\
\indent{\bf Theorem 3.1.~~}Suppose that  ($\mathbf{A}$) is valid and the ascent of $\lambda_{j,h_{m}}$ ($j=i,\cdots,i+q-1$) is 1 , and  there exist two eigenpairs ($\lambda_i,{u}_{i},{w}_{i}$) and ($\lambda^*_i,{u}^*_{i},{w}^*_{i}$) such that the eigenpair approximations $(\lambda_{j,h_{m}}$, ${u}_{j,h_{m}}$, ${w}_{j,h_{m}})$, $(\lambda^*_{j,h_{m}}$, ${u}^*_{j,h_{m}}$, ${w}^*_{j,h_{m}})\in \mathcal{C}\times \mathbf{H}_{h_{m}}$     have the following estimates
\begin{eqnarray}\label{s3.2}
  \| ({u}_{j,h_{m}},{w}_{j,h_{m}})- ({u}_j,w_j)\|_{\mathbf{H}}&\lesssim& \varepsilon_{h_{m}}(\lambda_i),\\\label{s3.3}
    \|({u}_{j,h_{m}},{w}_{j,h_{m}})- ({u}_j,w_j)\|_{\mathbf{H}_s}&\lesssim& H^{r_s}\varepsilon_{h_{m}}(\lambda_i),~s=0,1,\\\label{s3.4}
      \| ({u}^*_{j,h_{m}},{w}^*_{j,h_{m}})- ({u}^*_j,w^*_j)\|_{\mathbf{H}}&\lesssim& \varepsilon^*_{h_{m}}(\lambda_i),\\\label{s3.5}
    \|({u}^*_{j,h_{m}},{w}^*_{j,h_{m}})- ({u}^*_j,w^*_j)\|_{\mathbf{H}_s}&\lesssim& H^{r_s}\varepsilon^*_{h_{m}}(\lambda_i),~s=0,1,\\\label{s3.6}
    |\lambda_i-\lambda_{j,h_{m}}|&\lesssim& \varepsilon_{h_{m}}(\lambda_i)\varepsilon^*_{h_{m}}(\lambda_i).
\end{eqnarray}
Then after one correction step, there exists two eigenpairs $({\widehat{u}}_j,\widehat{w}_j)$, $({\widehat{u}}^*_j,\widehat{w}^*_j)$ such that the resultant approximations $(\lambda_{j,h_{m+1}},{u}_{j,h_{m+1}}, {w}_{j,h_{m+1}})$, $(\lambda^*_{j,h_{m+1}}$, ${u}^*_{j,h_{m+1}}$, ${w}^*_{j,h_{m+1}})\in \mathcal{C}\times \mathbf{H}_{h_{m+1}}$ has the following estimates
\begin{eqnarray}\label{s3.7}
 \| ({u}_{j,h_{m+1}},{w}_{j,h_{m+1}})- ({\widehat{u}}_j,\widehat{w}_j)\|_{\mathbf{H}}&\lesssim& \varepsilon_{h_{m+1}}(\lambda_i),\\\label{s3.8}
   \|({u}_{j,h_{m+1}},{w}_{j,h_{m+1}})- ({\widehat{u}}_j,\widehat{w}_j)\|_{\mathbf{H}_s}&\lesssim& H^{r_s}\varepsilon_{h_{m+1}}(\lambda_i),~s=0,1,\\\label{s3.9}
 \| ({u}^*_{j,h_{m+1}},{w}^*_{j,h_{m+1}})- ({\widehat{u}^*}_j,\widehat{w}^*_j)\|_{\mathbf{H}}&\lesssim& \varepsilon^*_{h_{m+1}}(\lambda_i),\\\label{s3.10}
   \|({u}^*_{j,h_{m+1}},{w}^*_{j,h_{m+1}})- ({\widehat{u}^*}_j,\widehat{w}^*_j)\|_{\mathbf{H}_s}&\lesssim& H^{r_s}\varepsilon^*_{h_{m+1}}(\lambda_i),~s=0,1,\\\label{s3.11}
    |\lambda_i-\lambda_{j,h_{m+1}}|&\lesssim& \varepsilon_{h_{m+1}}(\lambda_i)\varepsilon^*_{h_{m+1}}(\lambda_i),
\end{eqnarray}
where $$\varepsilon_{h_{m+1}}(\lambda_i):=\delta_{h_{m+1}}(\lambda_i)+H^{r_0}\varepsilon_{h_m}(\lambda_i)$$ and $$\varepsilon^*_{h_{m+1}}(\lambda^*_i):=\delta^*_{h_{m+1}}(\lambda^*_i)+H^{r_0}\varepsilon^*_{h_m}(\lambda^*_i).$$\\
\indent{\bf Proof.~~} Since $\Big\{({u}_{j,h_{m+1}},{w}_{j,h_{m+1}})\Big\}_{j=i}^{i+q-1}$ is a basis of $M_{h_m}(\lambda_i)$, $\Big\{({ {u}}_j, {w}_j)\Big\}_{j=i}^{i+q-1}$ is a basis of $M(\lambda_i)$.
For any $(v,z)\in M(\lambda_i)$ and $\|(v,z)\|_\mathbf{H}=1$ we denote
\begin{eqnarray*}
    (v,z)=\sum\limits_{j=i}^{i+q-1}\gamma_j ({ {u}}_j, {w}_j).
\end{eqnarray*}
For $(\widetilde{u}^*_{l,h_{m}},\widetilde{w}^*_{l,h_{m}})\in M_{h_m}(\lambda_i)$ in Assumption ($\mathbf{A}$), thanks to (\ref{s2.28}) there exists $(\widetilde{u}^*_l,\widetilde{w}^*_l)\in M(\lambda_i)$ satisfying
\begin{eqnarray*}
 \|(\widetilde{u}^*_l,\widetilde{w}^*_l)-(\widetilde{u}^*_{l,h_{m}},\widetilde{w}^*_{l,h_{m}}))\|_\mathbf{H}\lesssim\delta^*_{h_{m}}(\lambda^*_i).
\end{eqnarray*}

Then
\begin{eqnarray*}
  A((v,z),(\widetilde{u}^*_l,\widetilde{w}^*_l))=\sum\limits_{j=i}^{i+q-1}\gamma_j A(({ {u}}_j, {w}_j),(\widetilde{u}^*_l,\widetilde{w}^*_l)).
\end{eqnarray*}
Hence
\begin{eqnarray*}
|\gamma_l|&=&\frac{1}{A((u_l,w_l),(\widetilde{u}^*_l,\widetilde{w}^*_l))}\Big\{
  A((v,z),(\widetilde{u}^*_l,\widetilde{w}^*_l))-\sum\limits_{j=i,j\neq l}^{i+q-1} \gamma_j A(({ {u}}_j, {w}_j),(\widetilde{u}^*_l,\widetilde{w}^*_l))\Big\}.
\end{eqnarray*}
Due to Assumption $\mathbf{(A)}$ and $\|(\widetilde{u}^*_l,\widetilde{w}^*_l)-(\widetilde{u}^*_{l,h_{m}},\widetilde{w}^*_{l,h_{m}}))\|_\mathbf{H}\lesssim\delta^*_{h_{m}}(\lambda^*_i)$, since
\begin{eqnarray*}
|A(({ {u}}_j, {w}_j),(\widetilde{u}^*_l,\widetilde{w}^*_l))|&\le&|A(({ {u}}_j, {w}_j)-( {u}_{j,h_{m}}, {w}_{j,h_{m}}),(\widetilde{u}^*_l,\widetilde{w}^*_l))|\\
&&+|
  A(( {u}_{j,h_{m}}, {w}_{j,h_{m}}),(\widetilde{u}^*_l,\widetilde{w}^*_l)-(\widetilde{u}^*_{l,h_{m}},\widetilde{w}^*_{l,h_{m}}))|\\
  &&+|
  A(( {u}_{j,h_{m}}, {w}_{j,h_{m}}),(\widetilde{u}^*_{l,h_{m}},\widetilde{w}^*_{l,h_{m}}))|\\
  &\lesssim& \varepsilon_{h_{m}}(\lambda_i)+\delta^*_{h_{m}}(\lambda^*_i)+H^{r_0},
\end{eqnarray*}
we have
\begin{eqnarray*}
|\gamma_l|
  &\lesssim& \frac{1}{|A((u_l,w_l),(\widetilde{u}^*_l,\widetilde{w}^*_l))|}\Big\{
  1+\sum\limits_{j=i,j\neq l}^{i+q-1} |\gamma_j| (\varepsilon_{h_{m}}(\lambda_i)+\varepsilon^*_{h_{m}}(\lambda^*_i)+H^{r_0})\Big\}.
\end{eqnarray*}
Since $|A(({u}_{l,h_{m}},{w}_{l,h_{m}}),(\widetilde{u}^*_{l,h_{m}},\widetilde{w}^*_{l,h_{m}}))|$ has a positive lower bound uniformly with respect to $h_m$,
so does $|A(({u}_{l},{w}_{l}),({u}^*_{l},{w}^*_{l}))|$. It is immediate that
\begin{eqnarray*}
\sum\limits_{l=i}^{i+q-1}|\gamma_l|
  &\lesssim&
  q+q\sum\limits_{j=i}^{i+q-1} |\gamma_j| (\varepsilon_{h_{m}}(\lambda_i)+\varepsilon^*_{h_{m}}(\lambda^*_i)+H^{r_0}),
\end{eqnarray*}
from which it follows that
\begin{eqnarray*}
\sum\limits_{j=i}^{i+q-1}|\gamma_j|
  &\lesssim& 1.
\end{eqnarray*}
We set $\alpha_j=\lambda_j/\lambda_{j,h_m}$. By virtue of the orthogonality of $P_{h_{m+1}}$ and boundedness of $B(\cdot,\cdot)$ we have
\begin{eqnarray*}
&&\|\alpha_j(\widehat{u}_{j,h_{m+1}},\widehat{w}_{j,h_{m+1}})-P_{h_{m+1}}({u}_j,{w}_j)\|_\mathbf{H}^2\\
&&~~\lesssim A(\alpha_j(\widehat{u}_{j,h_{m+1}},\widehat{w}_{j,h_{m+1}})-P_{h_{m+1}}({u}_j,{w}_j),\alpha_j(\widehat{u}_{j,h_{m+1}},\widehat{w}_{j,h_{m+1}})-P_{h_{m+1}}({u}_j,{w}_j))\\
&&~~=A(\alpha_j(\widehat{u}_{j,h_{m+1}},\widehat{w}_{j,h_{m+1}})-({u}_j,{w}_j),\alpha_j(\widehat{u}_{j,h_{m+1}},\widehat{w}_{j,h_{m+1}})-P_{h_{m+1}}({u}_j,{w}_j))\\
&&~~=B(\lambda_{j,h_{m}}\alpha_j( {u}_{j,h_{m}}, {w}_{j,h_{m}})-\lambda_j({u}_j,{w}_j),\alpha_j(\widehat{u}_{j,h_{m+1}},\widehat{w}_{j,h_{m+1}})-P_{h_{m+1}}({u}_j,{w}_j))\\
&&~~\lesssim \|( {u}_{j,h_{m}}, {w}_{j,h_{m}})-({u}_j,{w}_j)\|_\mathbf{H_0}\|\alpha_j(\widehat{u}_{j,h_{m+1}},\widehat{w}_{j,h_{m+1}})-P_{h_{m+1}}({u}_j,{w}_j)\|_\mathbf{H},
\end{eqnarray*}
which together with  (\ref{s3.3}) yields
\begin{eqnarray*}
&&\|\alpha_j(\widehat{u}_{j,h_{m+1}},\widehat{w}_{j,h_{m+1}})-P_{h_{m+1}}({u}_j,{w}_j)\|_\mathbf{H}\lesssim H^{r_0}\varepsilon_{h_m}(\lambda_j).
\end{eqnarray*}
Thanks to Lemma 2.2, applying spectral approximation theory yields
\begin{eqnarray*}
 &&\| ({u}_{j,h_{m+1}},{w}_{j,h_{m+1}})- ({\widehat{u}}_j,\widehat{w}_j)\|_{\mathbf{H}} \lesssim \|({\widehat{u}}_j,\widehat{w}_j)-P_{h_{m+1}}({\widehat{u}}_j,\widehat{w}_j)\|_{\mathbf{H}}\\
 &&~~\lesssim
 \sup\limits_{(v,z)\in M(\lambda_i)\atop\|(v,z)\|_\mathbf{H}=1}\inf\limits_{(v_{h_{m+1}},z_{h_{m+1}})\in \mathbf{H}_{H,h_{m+1}}}\|(v,z)-(v_{h_{m+1}},z_{h_{m+1}})\|_\mathbf{H}\\
 &&~~\lesssim
 \sup\limits_{\gamma_j}\|\sum\limits_{j=i}^{i+q-1}\gamma_j (({ {u}}_j, {w}_j)-\alpha_j(\widehat{u}_{j,h_{m+1}},\widehat{w}_{j,h_{m+1}}))\|_\mathbf{H}\\
 &&~~\lesssim \sum\limits_{j=i}^{i+q-1}
\| ({ {u}}_j, {w}_j)-P_{h_{m+1}}({u}_j,{w}_j)+P_{h_{m+1}}({u}_j,{w}_j)-\alpha_j(\widehat{u}_{j,h_{m+1}},\widehat{w}_{j,h_{m+1}})\|_\mathbf{H}\\
   &&~~\lesssim
\delta_{h_{m+1}}(\lambda_j)+H^{r_0}\varepsilon_{h_m}(\lambda_j)=\varepsilon_{h_{m+1}}(\lambda_j),
 \end{eqnarray*}
 and for $s=0,1$,
 \begin{eqnarray*}
 &&\| ({u}_{j,h_{m+1}},{w}_{j,h_{m+1}})- ({\widehat{u}}_j,\widehat{w}_j)\|_{\mathbf{H_s}}\lesssim
  \|({\widehat{u}}_j,\widehat{w}_j)-P_{h_{m+1}}({\widehat{u}}_j,\widehat{w}_j)\|_{\mathbf{H}_s}\\
 &&~~\lesssim
  H^{r_s}\|({\widehat{u}}_j,\widehat{w}_j)-P_{h_{m+1}}({\widehat{u}}_j,\widehat{w}_j)\|_{\mathbf{H}}\\
  &&~~\lesssim
  H^{r_s} \varepsilon_{h_{m+1}}(\lambda_j),
 \end{eqnarray*}
 where we have used (\ref{2.20})  in the second inequality above.
The above argument implies (\ref{s3.7})-(\ref{s3.8}) hold. Similarly we can also prove (\ref{s3.9})-(\ref{s3.10}).  Using the proof method to show (\ref{s3.7}), we can derive (\ref{s3.11})  from (\ref{s2.30}).
~~~$\Box$\\

Now we are in a position to analyze the convergence of Multigrid Algorithm 2. \\
\indent{\bf Theorem 3.2.~~}Suppose   the conditions of Theorem 3.1 hold. Let the numerical eigenpair ($\lambda_{j,h_N}$, $u_{j,h_N}$, $w_{j,h_N},u^*_{j,h_N},w^*_{j,h_N})$ $(j=i,i+1,\cdots,i+q-1)$ be obtained by Algorithm 2. Then   there exists an eigenpair $(\lambda_i, {u_i},w_i, {u^*_i},w^*_i)$ such that   the following estimates hold
\begin{eqnarray}\label{s3.12}
  \| ({{u}}_{j,h_N},{w}_{j,h_N})-  ({{u}}_j,{w}_j)\|_{\mathbf{H}}&\lesssim& \delta_{h_{N}}(\lambda),\\\label{s3.13}
    \| ({{u}}_{j,h_N},{w}_{j,h_N})- ({{u}}_j,{w}_j)\|_{\mathbf{H}_s}&\lesssim& H^{r_s}\delta_{h_{N}}(\lambda), s=0,1,\\\label{s3.14}
      \| ({{u}}^*_{j,h_N},{w}^*_{j,h_N})-  ({{u}}^*_j,{w}^*_j)\|_{\mathbf{H}}&\lesssim& \delta^*_{h_{N}}(\lambda),\\\label{s3.15}
    \| ({{u}}^*_{j,h_N},{w}^*_{j,h_N})- ({{u}}^*_j,{w}^*_j)\|_{\mathbf{H}_s}&\lesssim& H^{r_s}\delta_{h_{N}}(\lambda), s=0,1,\\\label{s3.16}
    |\lambda_i-\lambda_{j,h_{N}}|&\lesssim& \delta_{h_{N}}(\lambda)\delta^*_{h_{N}}(\lambda).
\end{eqnarray}
\indent{\bf Proof.~~}
It's immediate from (\ref{s2.28})-(\ref{s2.30})  that
 \begin{eqnarray*}
  \|({{u}}_{j,h_1},{w}_{j,h_1})- ({{u}}_j,{w}_j)\|_\mathbf{H}&\lesssim& \delta_{h_{1}}(\lambda),\\
    \|({{u}}_{j,h_n},{w}_{j,h_1})- ({{u}}_j,{w}_j)\|_{\mathbf{H}_s}&\lesssim& H^{r_s}\delta_{h_{1}}(\lambda), s=0,1,\\
    |\lambda-\lambda_{h_{1}}|&\lesssim& \delta_{h_{1}}(\lambda)\delta^*_{h_{1}}(\lambda).
\end{eqnarray*}
Let $ \varepsilon_{h_{1}}(\lambda_i):= \delta_{h_{1}}(\lambda_i)$, due to $\delta_{h_{m}}(\lambda_i)\lesssim\beta\delta_{h_{m+1}}(\lambda_i)$, by recursion,
\begin{eqnarray*}
 \varepsilon_{h_{N}}(\lambda_i)&=&\delta_{h_{N}}(\lambda_i)+H^{r_0}\varepsilon_{h_{N-1}}(\lambda_i)\\
 &\lesssim&\delta_{h_{N}}(\lambda_i)+H^{2r_0}\delta_{h_{N-1}}(\lambda_i)+H^{r_0}\varepsilon_{h_{N-2}}(\lambda_i)\\
  &\lesssim&\sum\limits_{m=1}^{N} H^{(N-m)r_0}\delta_{h_{m}}(\lambda_i)\\
  &\lesssim& \sum\limits_{m=1}^{N} H^{(N-m)r_0}\beta^{N-m}\delta_{h_{N}}(\lambda_i)\\
   &\lesssim& \frac{1}{1-H^{r_0}\beta}\delta_{h_{N}}(\lambda_i)\\
   &\lesssim&  \delta_{h_{N}}(\lambda_i).
\end{eqnarray*}
Likewise, denoting $\varepsilon^*_{h_{1}}(\lambda_i):= \delta^*_{h_{1}}(\lambda_i)$ we can prove $\varepsilon^*_{h_{N}}(\lambda_i)  \lesssim \delta^*_{h_{N}}(\lambda_i)$.
   Using Theorem 3.1   we can obtain the desired results.\\

 \indent{\bf  Remark 3.1.}~~Similar to Section 5 of \cite{ji2}, we can  estimate the computational work
of Algorithm 2  and  prove
that solving the transmission eigenvalue problem needs almost the same
work as solving the corresponding boundary value problem.
 \section{Numerical experiment}

\indent In this section, we will report  some numerical
experiments for solving the transmission eigenvalue problem (\ref{s2.6}) by
multigrid method (Algorithm 2)
to validate our theoretical results.\\
\indent We use Matlab 2012a to solve (\ref{s2.1})-(\ref{s2.4}) on a
Lenovo G480 PC with 4G memory. Our program is compiled under the
package of iFEM \cite{chen}. In Algorithms 2, we use internal command $eigs$ in Matlab to solve matrix eigenvalue problem  use command $'\backslash'$ in
Matlab to solve the associated linear algebraic systems. In numerical examples, since the exact eigenvalues are   unknown, 
we compute the relatively accurate  ones to replace them.\\
\indent In our numerical experiments, we construct the finite element space
$\mathbf{H}_{H,h_{m+1}}$ as $U\times V=:$ $$\Big(S^H+ span\Big\{\widehat{u}_{j,h_{m+1}}, \widehat{u}^*_{j,h_{m+1}}\Big\}_{j=i}^{i+q-1}\Big)\times \Big(S^H+ span\Big\{\widehat{w}_{j,h_{m+1}},\widehat{w}^*_{j,h_{m+1}}\Big\}_{j=i}^{i+q-1}\Big),$$
so that the first equation in Step 2 of Algorithm 1 can be conveniently written as
\begin{eqnarray*}
&&(\frac{1}{n-1}\Delta u_{j,h_{m+1}}, \Delta v)_{0}=\lambda_{j,h_{m+1}}\Big\{
(\nabla(\frac{1}{n-1}u_{j,h_{m+1}}),
\nabla v)_{0}\nonumber\\
&&~~~~~~+(\nabla {u}_{j,h_{m+1}}, \nabla(\frac{n}{n-1}{v} ))_{0}-
(\frac{n}{n-1}w_{j,h_{m+1}}, v  )_{0}\Big\},~~\forall v
\in~U,\\
 &&(w_{j,h_{m+1}},
z)_{0}=\lambda_{j,h_{m+1}}(u_{j,h_{m+1}},z )_{0},~~\forall z \in V.
\end{eqnarray*}
%
\indent For reading conveniently, we use the following notations in our tables and figures:\\
\indent $k_{j,h}=\sqrt{\lambda_{j,h}}$: The $jth$ eigenvalue
obtained by
  Algorithm 2 on $\mathbf{H}_{H,h}$.\\


\indent We  consider the case when $D$ is the unit square $(0,1)^2$ or the L-shaped domain $(-1,1)^2\backslash [0,1)\times (-1,
0]$ and
the index of refraction $n=16,4,8+x_{1}-x_{2}$. We choose $S^h$ to be the
Bogner-Fox-Schmit rectangle element (BFS element) space such that $S^h\subset H^2_0(D)$, and the numerical results are shown in Tables 1-2.
Besides, we depict the error curves for the numerical eigenvalues
(see Figures 1-3) based on the
results in Tables 1-2.\\
\indent According to the regularity theory,  we know $u, w\in
H^{4}(D)$ if $D$ is a square. When the ascent of $k$ is equal to $1$: according to (\ref{s3.16}), the convergence order of the eigenvalue
approximation $k_{j,h}$ is four.  It is seen from  Figures 1-3 that the convergence order of the
numerical eigenvalues on the unit square domain is four, which
coincides with the theoretical result. \\
\indent Figures 1-3 indicates that  on the L-shaped domain, the
convergence order of the listed numerical eigenvalues is less than four.
  This fact suggests that their corresponding  eigenfunctions on the non-convex domain do
have singularities in different degrees.\\
\indent In addition, Tables 1-2 show that the numerical
eigenvalues obtained by Algorithm 2 using BFS elements give a good approximation; our multigrid method can achieve   the
same convergence order as the one in \cite{ji2}. 
It is
worth  noticing that  our method can   compute real and complex transmission eigenvalues efficiently (see Tables 1-2). 
\\\\
\indent${Submitted~to~Applied~Mathematics~and~Computation ~on~Oct~17,~2015}$

\begin{table}
\caption{The numerical eigenvalues obtained by Algorithm 2 on the unit square
domain with $H=h_1=\frac{\sqrt2}{8}$($n=16,8+x_1-x_2$) and $H=h_1=\frac{\sqrt2}{16}$($n=4$).}
\begin{center} \footnotesize
\begin{tabular}{llllllll}\hline
$n$&$h$&$k_{1,h}$&$k_{2,h},k_{3,h}$&$k_{4,h}$\\\hline
16&  $\frac{\sqrt2}{8}$&    1.880051827&	2.446255515&	2.868193148\\
16&  $\frac{\sqrt2}{16}$&    1.879621633&	2.444371226&	2.866560541\\
16&  $\frac{\sqrt2}{32}$&   1.879593109&	2.444244719&	2.866446979\\
16&  $\frac{\sqrt2}{64}$&   1.879591295&	2.444236640&	2.866439605\\
16&  $\frac{\sqrt2}{128}$&  1.879591180&	2.444236133&	2.866439141\\
16&  $\frac{\sqrt2}{256}$&   1.879591166&	2.444236099&	2.866439110
\\\hline\hline
$n$&$h$&\multicolumn{2}{l}{$k_{1,h}$,$k_{2,h}$}&$k_{3,h}$&$k_{5,h}$\\\hline
4&$\frac{\sqrt2}{16}$&\multicolumn{2}{l}{4.271570823 $\pm$ 1.147502410i}& 	    5.477120630&	6.100857372\\
4&$\frac{\sqrt2}{32}$&\multicolumn{2}{l}{4.271689022 $\pm$ 1.147437637i}&	 	5.476172735&	6.100321070\\
4&$\frac{\sqrt2}{64}$&\multicolumn{2}{l}{4.271696373 $\pm$ 1.147433642i}&	 	5.476112619&	6.100286894\\
4&$\frac{\sqrt2}{128}$&\multicolumn{2}{l}{4.271696834 $\pm$ 1.147433395i}&	 	5.476108843&	6.100284742\\
4&$\frac{\sqrt2}{256}$&\multicolumn{2}{l}{4.271696869 $\pm$ 1.147433423i}&	 	5.476108632&	6.100284617\\
\hline\hline
$n$&$h$&$k_{1,h}$&$k_{2,h}$&$k_{5,h},k_{6,h}$\\\hline
$8+x_1-x_2$&  $\frac{\sqrt2}{8}$&    2.823445794&	3.542452244&	  \multicolumn{2}{l}{4.4971031374 $\pm$ 0.8770188489i}\\
$8+x_1-x_2$&  $\frac{\sqrt2}{16}$&   2.822270903&	3.538946571&	  \multicolumn{2}{l}{4.4966538580 $\pm$ 0.8718551076i}\\
$8+x_1-x_2$&  $\frac{\sqrt2}{32}$&   2.822194508&	3.538712669&	  \multicolumn{2}{l}{4.4965518150 $\pm$ 0.8714987351i}\\
$8+x_1-x_2$&  $\frac{\sqrt2}{64}$&   2.822189665&	3.538697695&	  \multicolumn{2}{l}{4.4965525625 $\pm$ 0.8714831407i}\\
$8+x_1-x_2$&  $\frac{\sqrt2}{128}$&  2.822189362&	3.538696758&	  \multicolumn{2}{l}{4.4965519527 $\pm$ 0.8714818531i}\\
$8+x_1-x_2$&  $\frac{\sqrt2}{256}$&   2.822189348&	3.538696701&	  \multicolumn{2}{l}{4.4965519517 $\pm$ 0.8714817941i}\\
\hline
\end{tabular}
\end{center}
\end{table}

\begin{table}
\caption{The numerical eigenvalues obtained by Algorithm 2 on the L-shaped
domain with $H=h_1=\frac{\sqrt2}{8}$.}
\begin{center} \footnotesize
\begin{tabular}{llllll}\hline
$n$&$h$&$k_{1,h}$&$k_{2,h}$&$k_{3,h}$&$k_{4,h}$\\\hline
16&  $\frac{\sqrt2}{8}$&    1.4850654&	1.5705634&	1.7078129&	1.7834681\\
16&  $\frac{\sqrt2}{16}$&   1.4802424&	1.5699011&	1.7061982&	1.7831490\\
16&  $\frac{\sqrt2}{32}$&   1.4780404&	1.5697716&	1.7055795&	1.7831209\\
16&  $\frac{\sqrt2}{64}$&   1.4770116&	1.5697385&	1.7052950&	1.7831171\\
16&  $\frac{\sqrt2}{128}$&  1.4765288&	1.5697294&	1.7051613&	1.7831163\\
\hline\hline
$n$&$h$&\multicolumn{2}{l}{$k_{1,h}$,$k_{2,h}$}&$k_{3,h}$&$k_{4,h}$\\\hline
4&$\frac{\sqrt2}{8}$&\multicolumn{2}{l}{3.1106061 $\pm$ 1.1986040i}&	 3.5307679&	3.6418388\\
4&$\frac{\sqrt2}{16}$&\multicolumn{2}{l}{3.1061706 $\pm$ 1.1925668i}&	 3.5239697&	3.6391912\\
4&$\frac{\sqrt2}{32}$&\multicolumn{2}{l}{3.1038283 $\pm$ 1.1897655i}&	 3.5213679&	3.6387478\\
4&$\frac{\sqrt2}{64}$&\multicolumn{2}{l}{3.1027011 $\pm$ 1.1884482i}&	 3.5202067&	3.6386426\\
4&$\frac{\sqrt2}{128}$&\multicolumn{2}{l}{3.1021669 $\pm$ 1.1878286i}&	 3.5196692&	3.6386141\\
\hline\hline
$n$&$h$&$k_{1,h}$&$k_{2,h}$& {$k_{5,h},k_{6,h}$}\\\hline
$8+x_1-x_2$&  $\frac{\sqrt2}{8}$&    2.3127184&	2.3974892&	\multicolumn{2}{l}{2.9287347 $\pm$ 0.5743458i}\\
$8+x_1-x_2$&  $\frac{\sqrt2}{16}$&   2.3069617&	2.3960567&	\multicolumn{2}{l}{2.9272395 $\pm$ 0.5686048i}\\
$8+x_1-x_2$&  $\frac{\sqrt2}{32}$&   2.3043809&	2.3957864&	\multicolumn{2}{l}{2.9254909 $\pm$ 0.5666799i}\\
$8+x_1-x_2$&  $\frac{\sqrt2}{64}$&   2.3031831&	2.3957182&	\multicolumn{2}{l}{2.9247531 $\pm$ 0.5656800i}\\
$8+x_1-x_2$&  $\frac{\sqrt2}{128}$&  2.3026210&	2.3956994&	\multicolumn{2}{l}{2.9245398 $\pm$ 0.5649109i}\\

\hline
\end{tabular}
\end{center}
\end{table}

\newpage

  \begin{figure}
\includegraphics[width=0.45\textwidth]{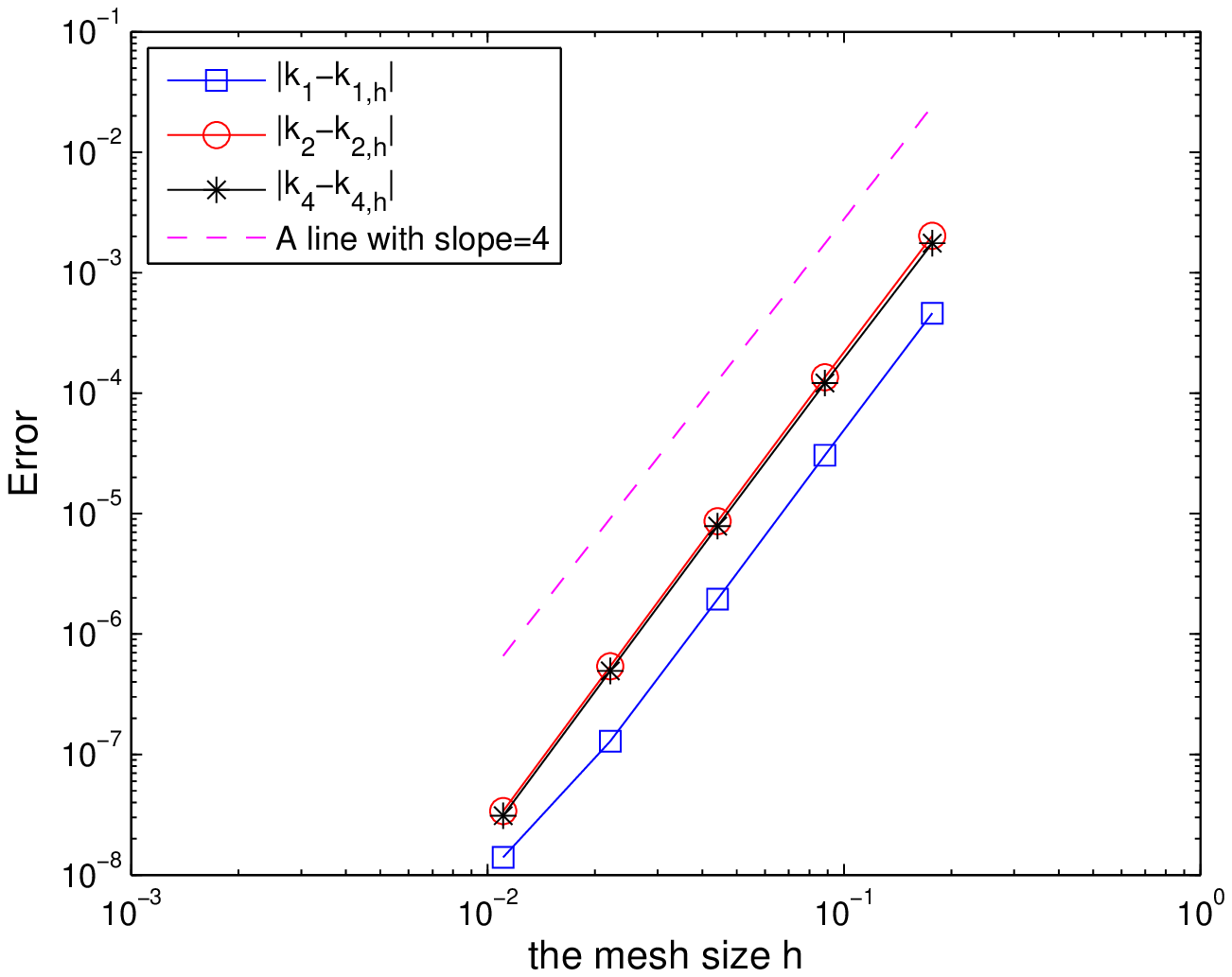}
 \includegraphics[width=0.45\textwidth]{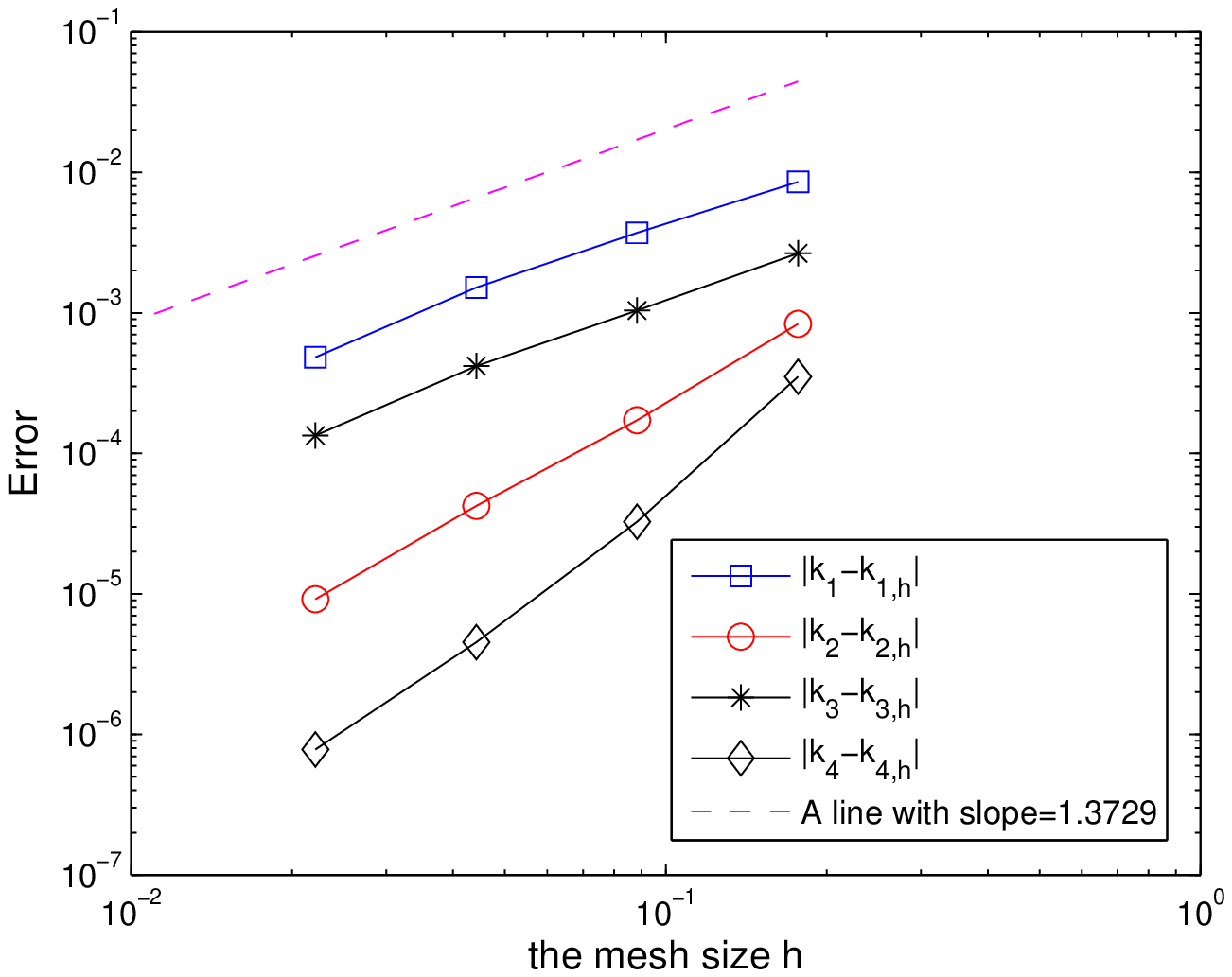}
\caption{{ Error curves for numerical eigenvalues on the unit square(left) and on the L-shaped(right) with
$n=16$.}}
 \end{figure}
  \begin{figure}
\includegraphics[width=0.45\textwidth]{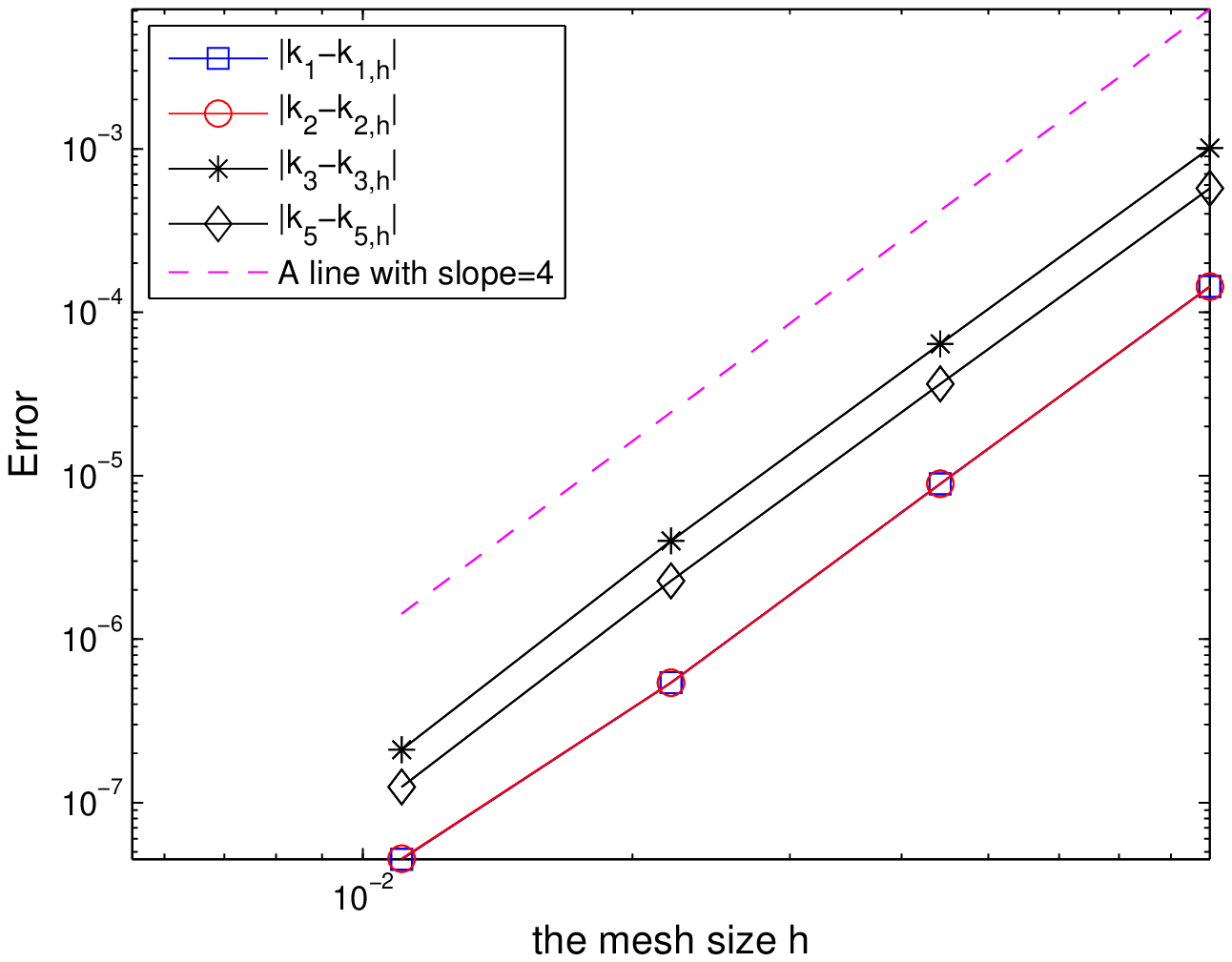}
 \includegraphics[width=0.45\textwidth]{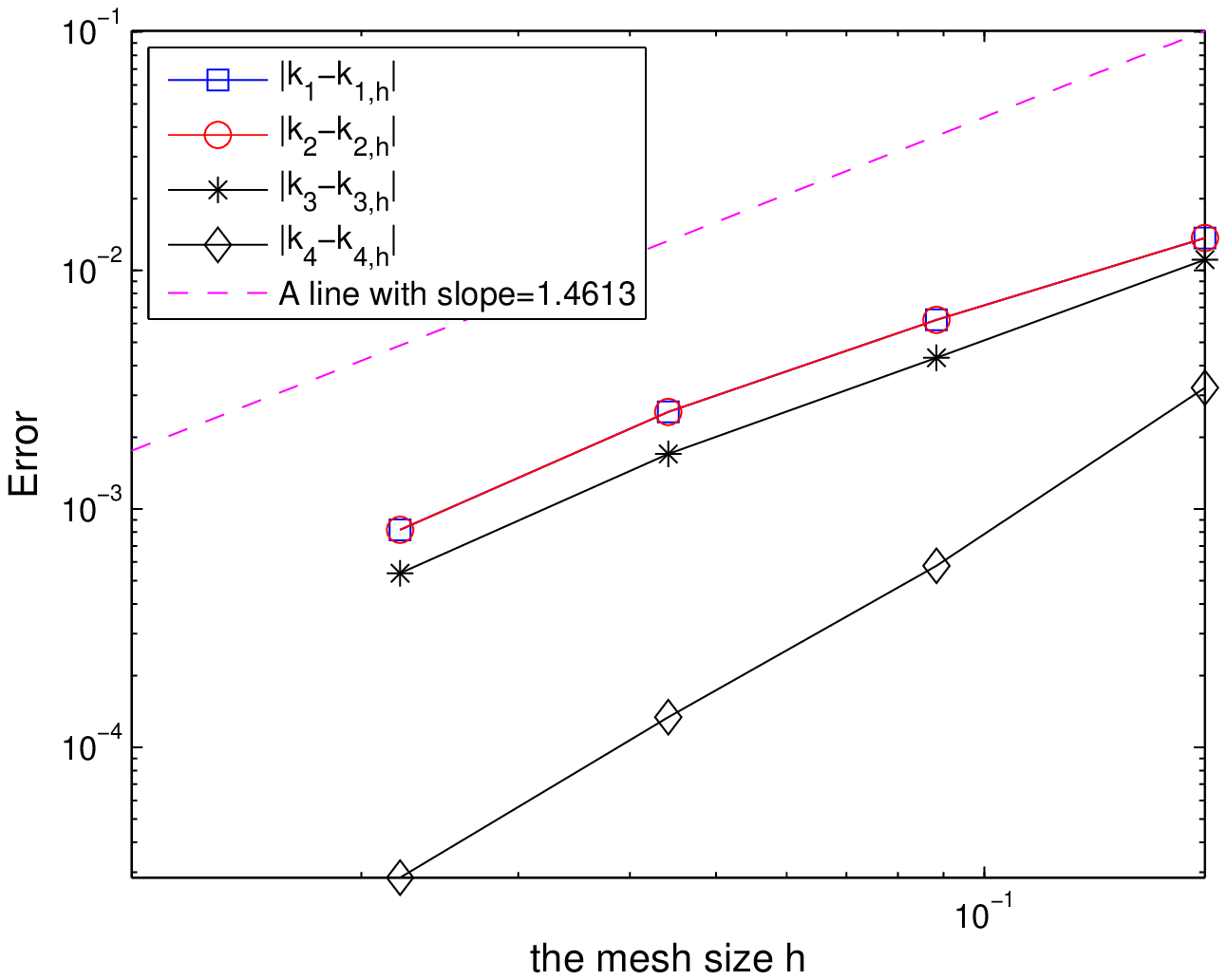}
\caption{{ Error curves for numerical eigenvalues on the unit square(left) and on the L-shaped(right) with
$n=4$.}}
 \end{figure}
  \begin{figure}
  \includegraphics[width=0.45\textwidth]{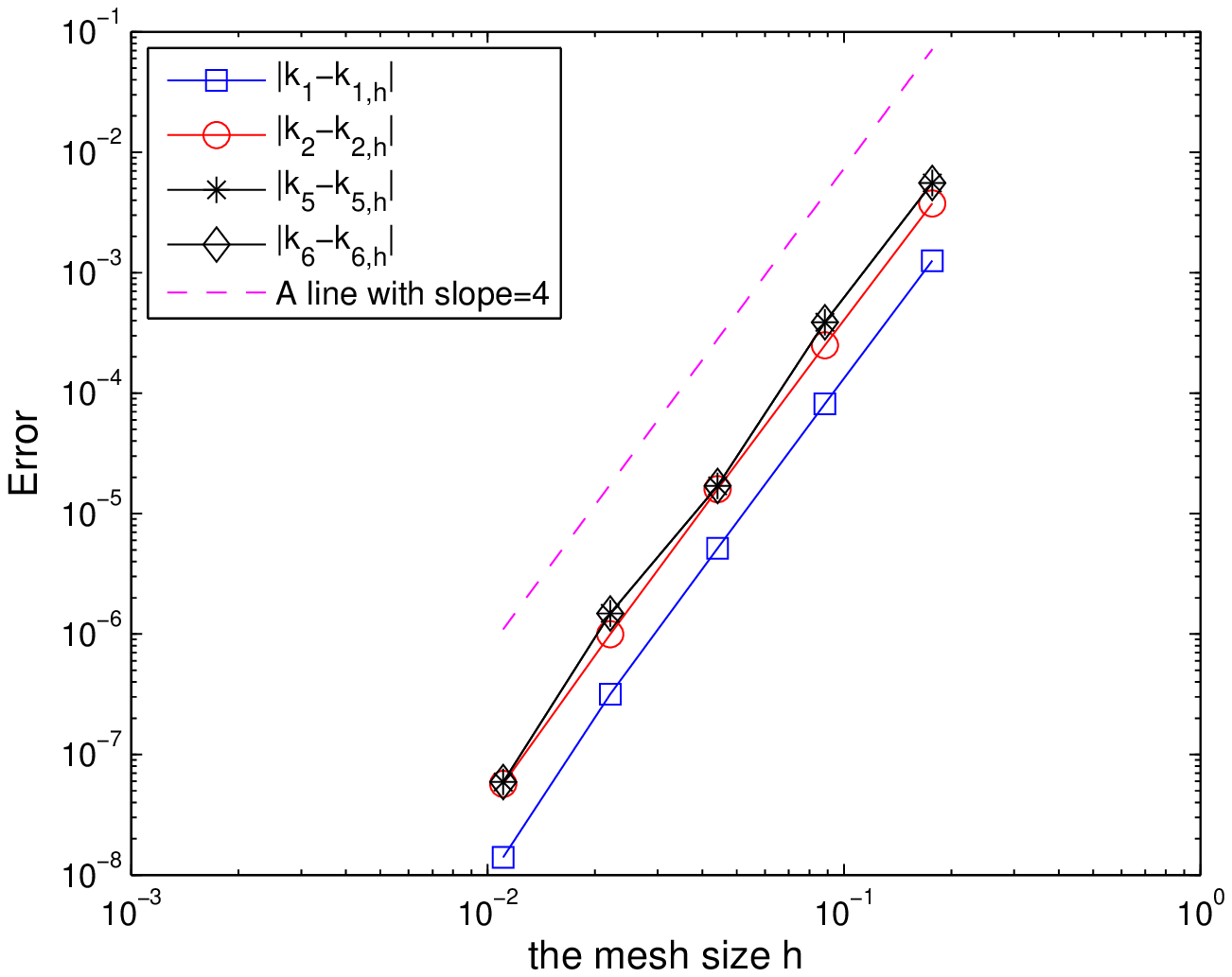}
 \includegraphics[width=0.45\textwidth]{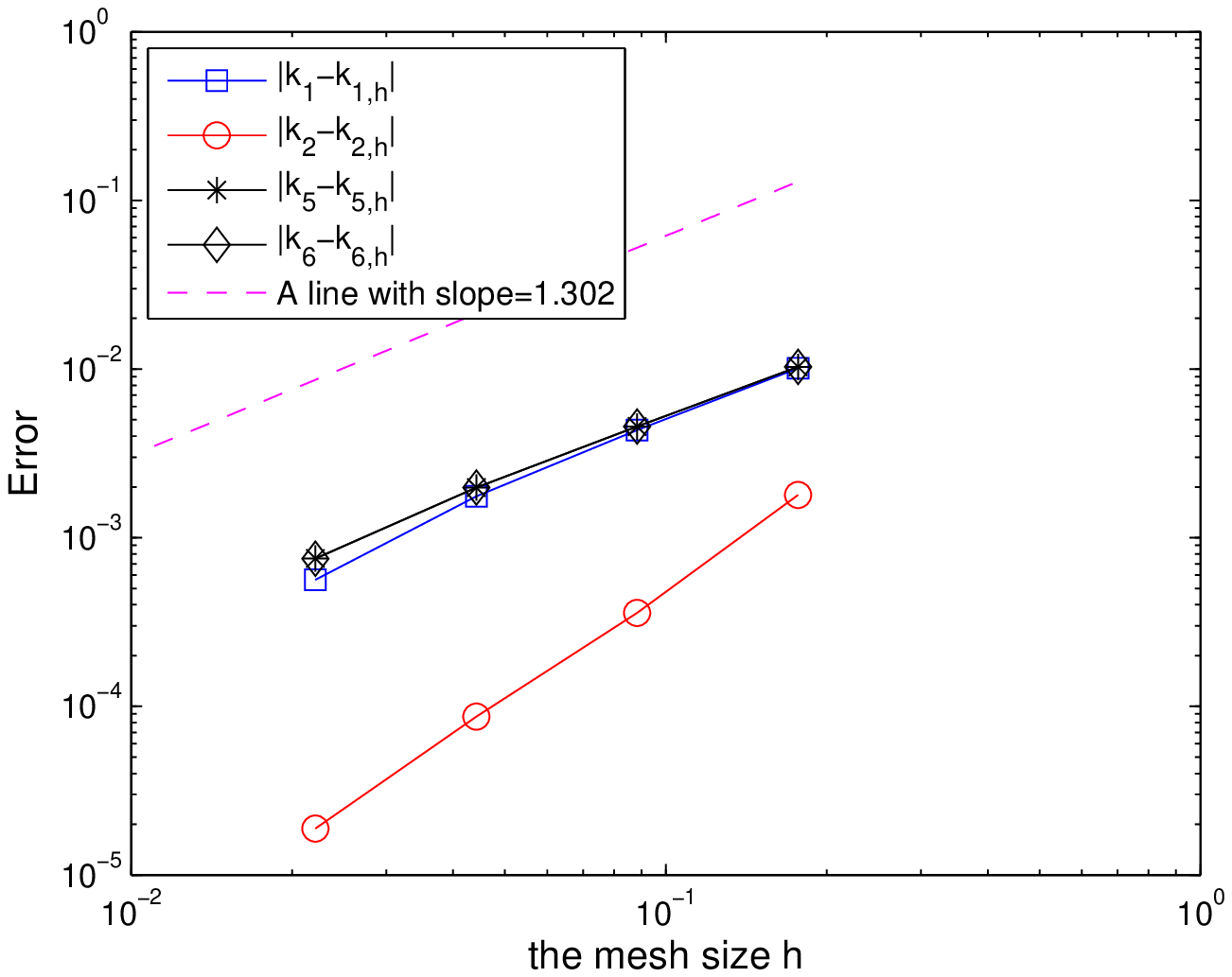}
 \caption{{ Error curves for numerical eigenvalues on the unit square(left) and on the L-shaped(right) with
$n=8+x_1-x_2$.}}
  \end{figure}
\end{document}